\documentclass{article}
\usepackage[utf8]{inputenc}
 \usepackage{esvect} 
 \usepackage{latexsym}
\usepackage{amsfonts}       
\usepackage{color}
\usepackage{amsmath}
\newcounter{eqnsave}

\def\C{{{{\rm {\mbox{\small l}}} \kern -.50em {\rm C}}}}
\def\I{{{{\rm l} \kern -.10em {\rm I}}}}
\def\R{{{{\rm l} \kern -.15em {\rm R}}}}
\def\N{{{{\rm l} \kern -.15em {\rm N}}}}
\def\E{{{{\rm l} \kern -.15em {\rm E}}}}

\newcommand{\bea}{\begin{eqnarray}}
\newcommand{\beas}{\begin{eqnarray*}}

\newcommand{\ba}{\begin{array}}
\newcommand{\ea}{\end{array}}

\newtheorem{theorem}{Theorem}[section]
\newtheorem{lemma}[theorem]{Lemma}

\newtheorem{definition}{Definition}[section]

\newcommand{\qed}{{}\hfill $\Box$ \\[15pt]}
\newcommand{\eea}{\end{eqnarray}}                %
\newcommand{\bas}{\begin{eqnarray*}}            %
\newcommand{\eas}{\end{eqnarray*}}              %
\title{On the existence of solutions for a parabolic-elliptic   chemotaxis model with flux limitation and logistic source}

\author{S. Sastre-Gomez\footnote{
 Departamento de  Ecuaciones Diferenciales y An\'alisis Num\'erico,  Facultad de Matem\'aticas, Universidad  de Sevilla. e-mail: ssastre@us.es}  \ and   J. Ignacio Tello\footnote{
 Departamento de Matem\'aticas Fundamentales,  Facultad de Ciencias, Universidad Nacional de Educaci\'on a Distancia. e-mail: jtello@mat.uned.es}}

\begin{document}

\maketitle

\begin{abstract}
In this article we study the existence of   solutions of  a parabolic-elliptic system of partial differential equations  
describing the behaviour of a biological species ``$u$" and a chemical stimulus
 ``$v$"
 in a bounded and regular domain $\Omega$ of $\R^N$.  
  The equation for $u$ is a parabolic equation with a nonlinear second order term of chemotaxis type 
 with  flux limitation   as  
  $$ -\chi div (u  |\nabla  \psi|^{p-2} \nabla v), $$ 
  for $p>1$. 
  The chemical substance distribution  $v$  satisfies the elliptic equation 
  $$-\Delta v+v=u.$$
  The evolution of $u$ is also determined by a logistic type growth term $\mu u(1-u)$. The system is studied under  homogeneous Neumann boundary conditions.  
  The main result of the article is the existence of uniformly bounded solutions for  $p<3/2$ and any  $N\ge 2$.
\end{abstract}

\section{Introduction}

In the decade of 70$^{\prime}$s  of the last century,  Keller and Segel  \cite{ks} and \cite{ks2} described  chemotaxis phenomena for a biological species density ``$u$" and a chemical substance concentration ``$v$" using systems of partial differential equations. 
The Keller-Segel model presents a second order chemotaxis term with a linear dependence  of the chemical flux ``$\nabla v$" in the way 
$$- div (\chi u \nabla v)$$ 
where $\chi$ is a given constant.  
After the model was proposed, many authors have studied such systems from a  mathematical point of view. 
See for instance
  Horstmann \cite{horstmann1}, \cite{horstmann2} and Hillen and Painter \cite{hp}   for details concerning  mathematical results of chemotaxis models with linear flux. 

 In  M.A. Rivero, R.T. Tranquillo, H.M. Buettner and  D.A. Lauffenburger \cite{1989}, the authors proposed a   chemotaxis term with non-linear dependence of the chemical flux, given in the form 
$$- \frac{\partial }{\partial x}  \left(   F(u, \frac{\partial v}{\partial x}  )\right), \qquad  F(u, \frac{\partial v}{\partial x}  )=
 arctan \left[ \frac{u}{1+u} \frac{\partial v}{\partial x} \right]
$$
in a  one-dimensional spatial domain. 

   A. Chertock, A. Kurganov, X. Wang and Y. Wu  in   \cite{chertov} consider the  
parabolic-parabolic system with a flux limitation defined by
$$F(  \nabla v  )=
  \frac{\nabla v}{\sqrt{1+ |\nabla v|^2}}.
$$

 A. Bianchi,  K.J. Painter and  J.A. Sherratt \cite{bps1} and \cite{bps2}  use such a  nonlinear dependence to obtain a biomedical model.   
 N. Bellomo and  M. Winkler in \cite{bw1} and \cite{bw2} study the parabolic-elliptic system  for  a chemotaxis term of the form 
$$- div (u f(|\nabla v|) \nabla v )$$ 
for 
$$f(|\nabla v |)=\frac{\chi}{\sqrt{1+|\nabla v|^2}} $$ 
and a  diffusive term  
$$-div\left( \frac{u \nabla u}{\sqrt{u^2+|\nabla u|^2}}\right)$$
to describe the evolution of the biological species $u$. 

In \cite{bw1}, for $\chi <1$ and $N \geq 2$ the solutions are global and  bounded  if the initial mass is smaller than $(\chi^2-1)_+^{-\frac{1}{2}}$.  Using an energy method,   Winkler \cite{winkler} has  obtained blow up of solutions for 
 $$f=f(|\nabla v |)=\frac{\chi}{(1+|\nabla v|^2)^{\frac{\alpha}{2}}}$$  when 
$$  \alpha<\frac{N-2}{2(N-1)}. $$

In this article we consider  a bounded      $N$-dimensional domain,    $\Omega $,   with regular boundary $\partial \Omega$ 
and
denote by $\vv{n}$ the  outward pointing normal vector on the boundary $\partial \Omega$.

The equation for $v$ is restricted to the elliptic case, for  simplicity, we assume that $v$ satisfies the Poisson equation and  the system studied reads as  follows
\begin{align}{} 
\label{1.1} 
 &u_t-\Delta u = -  div  (\chi u|\nabla v|^{p-2}\nabla v)+\mu u(1-u)  & x\in\Omega,\quad t>0, &
 \\[2mm]  \label{1.2}
 & -\Delta v  +v= u &  x\in\Omega,\quad t>0,  &
  \\[2mm]  \label{1.3}  &
 \displaystyle    \frac{\partial  u}{\partial \vv{n} }=\displaystyle\frac{\partial v}{\partial \vv{n} }=0 &  x\in \partial \Omega,\quad t>0,  &
\\[2mm]   & u(0,x)=u_0(x),~~~~ v(0,x)=v_0(x)  &   x\in\Omega.  \label{1.4} &
     \end{align}      
In M. Negreanu and J.I. Tello \cite{nt8}, the system is considered for $\mu=0$, i.e., the logistic term does not appear and  $p$ satisfies
$$p \in (1, \infty), \quad \mbox{ if } N=1 \quad \mbox{ and } \quad p\in \left(1, \frac{N}{N-1}\right),   \quad \mbox{ if } N\geq 2,$$
to   obtain uniform bounds in $L^{\infty}(\Omega)$.

 The  complementary case,
  for $p\in (N/(N-1),2),$ for $N>2$ and  $\Omega$ defined as  the unit ball, presents blow up for some initial data,  see \cite{t}. 

In this article we study  the global existence of solutions  in $(0,T)$ for any $T<\infty$ under the 
following assumptions   
\begin{equation} \label{H1} 
\Omega \mbox{ is an open and bounded domain with regular boundary }  \partial \Omega, 
\end{equation} 
\begin{equation} 
\label{H2}
\left\{ 
\begin{array}{ll} 
p<2, & N=2,
\\ [2mm] 
p \in (1, 3/2), & N \geq 3,
\end{array}
\right.
\end{equation} 
and the initial data $u_0$  satisfy 
\begin{equation}  \label{H3} 
u_0 \in C^{2+ \alpha} (\Omega), \quad  
   \frac{\partial  u_0}{\partial \vv{n} } =0, \quad    x\in \partial \Omega,
    \end{equation}
   for some  $\alpha \in (0,1)$. 
    
    In this article we prove the existence of weak solutions to problem (\ref{1.1})-(\ref{1.4}) in the  sense of the following definition.
\begin{definition}\label{def1.1}
We say that $(u,v)\in L^2(0,T:H^1(\Omega))^2$ is a weak solution to \eqref{2.1}-\eqref{2.4} if 
$u\nabla v \in [L^{1}(\Omega)]^N $
and 
for any $\varphi,\psi\in L^2(0,T:H^1(\Omega))\cap H^1(0,T:L^2(\Omega))$, 
$(u,v)$ satisfies 
 \begin{equation}
\begin{array}{r} 
\displaystyle \int_\Omega \varphi(t,x) u (t,x)dx- \int_0^t\int_\Omega 
u(s,x)\varphi_t(s,x)dxds+
 \int_0^t \int_\Omega \nabla u(s,x)\nabla \varphi(s,x)dxds=
 \\ [4mm] 
\displaystyle \int_\Omega u_0(x)\varphi(0,x)dx+\int_0^t\int_\Omega \chi u(s,x) |\nabla v(s,x)|^{p-1}
\nabla v_n(s,x) dx ds
\\
\displaystyle+\int_0^t\int_\Omega \mu u(s,x) (1-(u(s,x))_+)
\varphi(s,x)dx ds, 
     \end{array}
 \end{equation}
 \begin{equation}
\int_0^t\int_\Omega \nabla v(s,x)\nabla\psi(s,x)dx ds =
\int_0^t\int_\Omega \left( u(s,x)- v(s,x)\right)\psi(s,x) dx ds.
 \end{equation}
\end{definition} 
    The   existence of weak solutions given in the previous definition is  enclosed in the following theorem. 
    \begin{theorem} \label{t1} 
    Under assumptions (\ref{H1})-(\ref{H3}), there exists at least a solution to (\ref{1.1}) in $(0,T)$ for any $T<\infty$.  
    \end{theorem} 
    The article is organized as follows. In Section \ref{s2}, we obtain some a priori estimates that are used in Section \ref{s3} to proof Theorem \ref{t1}. We use an iterative method based in the Moser-Alikakos  iteration which allows us to obtain explicit  $L^q$ estimates of the solutions. Then,  we pass to the limit to get the   boundedness of the solutions in  $L^{\infty}$ norm. Finally, by using an approximated problem we pass to the limit in the weak formulation to obtain the existence of weak solutions.  
    
\section{A priori estimates for the auxiliary problem}  
  \label{s2} 
    \setcounter{equation}{0}

We  introduce the following auxiliary problem 
\begin{align}{} 
\label{2.1} 
 &u_{nt}-\Delta u_n = -  div  \left(\chi u_n\frac{|\nabla v_n|^{p-2}\nabla v_n}{1+\frac{1}{n} |\nabla  v_n|^{p-1}}\right)+\mu u_n(1-u_n)  & x\in\Omega,\quad t>0, &
 \\[2mm]  \label{2.2}
 & -\Delta v_n  +v_n= u_n &  x\in\Omega,\quad t>0,  &
  \\[2mm]  \label{2.3}  &
 \displaystyle    \frac{\partial  u_n}{\partial \vv{n} }=\displaystyle\frac{\partial v_n}{\partial \vv{n} }=0 &  x\in \partial \Omega,\quad t>0,  &
\\[2mm]   & u_n(0,x)=u_0(x), \   \   v_n(0,x)=v_0(x)  &   x\in\Omega.  \label{2.4} &
     \end{align}      %

\begin{definition}\label{def3.1}
We say that $(u_n,v_n)\in L^2(0,T:H^1(\Omega))^2$ is a weak solution to \eqref{2.1}-\eqref{2.4} if for any $\varphi,\psi\in L^2(0,T:H^1(\Omega))\cap H^1(0,T:L^2(\Omega))$, 
$(u_n,v_n)$ satisfies 
 \begin{equation}
\begin{array}{r} 
\displaystyle \int_\Omega \varphi(t,x) u_n (t,x)dx- \int_0^t\int_\Omega 
u_n(s,x)\varphi_t(s,x)dxds+
 \int_0^t \int_\Omega \nabla u_n(s,x)\nabla \varphi(s,x)dxds=
 \\ [4mm] 
\displaystyle \int_\Omega u_0(x)\varphi(0,x)dx+\int_0^t\int_\Omega \chi u_n(s,x)\frac{|\nabla v_n(s,x)|^{p-1}
\nabla v_n(s,x)}{1+\frac{1}{n}|\nabla v_n(s,x)|} \nabla\varphi(s,x)dx ds
\\
\displaystyle+\int_0^t\int_\Omega \mu u_n(s,x) (1-(u_n(s,x))_+)
\varphi(s,x)dx ds, 
     \end{array}
 \end{equation}
 \begin{equation}
\int_0^t\int_\Omega \nabla v_n(s,x)\nabla\psi(s,x)dx ds =
\int_0^t\int_\Omega \left( u_n(s,x)- v_n(s,x)\right)\psi(s,x) dx ds.
 \end{equation}

\end{definition}

In this section  we give some estimates that will be useful to prove the existence of solution of \eqref{2.1}-\eqref{2.4}. In particular we obtain uniform bounds in $L^\infty(\Omega)$ for $u_n$ and $v_n$. 

In the following result we obtain an upper bound of the total mass of $u_n$. 
\begin{lemma}  \label{l3.1}
The total mass of the component $u_n$ of the solution to 
(\ref{2.1}) is bounded,  $\displaystyle 
\int_{\Omega}u_n\, dx \leq c_1$, with $0< c_1 <\infty$.
\end{lemma}
{\bf Proof.} 
After integration by parts in equation (\ref{2.1}) we get
$$ \frac{d}{dt}\int_{\Omega} u_n\, dx =
\mu \int_{\Omega} u_n(1-u_n)dx .$$
Thanks to H\"older inequality we obtain 
$$\int_{\Omega} u_n^2 dx \geq |\Omega|^{-1}  \left(\int_{\Omega}
|u_n| ~dx \right)^2,
$$
and therefore
$$ \frac{d}{dt}\int_{\Omega} u_n\, dx \leq 
\mu \int_{\Omega} u_n\, dx - |\Omega|^{-1} \mu \left(\int_{\Omega}  |u_n|~dx\right)^2
\leq 
\mu \int_{\Omega} u_n\, dx - |\Omega|^{-1} \mu \left(\int_{\Omega}  u_n~dx\right)^2 
.$$
Thanks to Maximum principle for ODEs,  we have the result for 
$$c_1=:\max\left\{ \int_{\Omega} u_0~dx,~ |\Omega|\right\}.
$$
$\hfill\square$
\vspace{0.3cm}

 
 In the following result we prove that $u_n$ is non negative if the initial data $u_0$  is non negative. 

\begin{lemma}\label{l3.1.1}
 The solution  $u_n$ to 
(\ref{2.1})-\eqref{2.4} with $u_0\ge 0$  satisfies  $\displaystyle 
 u_n \geq 0$. 
\end{lemma}
\begin{proof}  We first consider the auxiliary equation 
\begin{equation}\label{2.5a} u_{nt}-\Delta u_n = -  div  \left(\chi u_n\frac{|\nabla v_n|^{p-2}\nabla v_n}{1+\frac{1}{n} |\nabla  v_n|^{p-1}}\right)+\mu u_n(1-(u_n)_+) 
\end{equation}
 and define $T_h$   as follows
 $$
 T_h(s)=\left\{
 \begin{array}{rll}
 -h, &\mbox{ if } & s\le -h, \\
 s, &\mbox{ if } &-h<s<0,\\
 0 ,& \mbox{ if } &s\ge 0.
 \end{array}\right.
 $$
 We denote the 
  primitive of $T_h$, by $\Phi_h$, which is given by
  $$
 \Phi_h(s)=\left\{
 \begin{array}{rll}
 -hs - \frac{h^2}{2} , &\mbox{ if } & s\le -h, \smallskip \\
 [2mm] 
 \frac{s^2}{2}, &\mbox{ if } &-h<s<0,\\ [2mm] 
 0 , & \mbox{ if } &s\ge 0,
 \end{array}\right.
 $$
and satisfies $\Phi_h^{\prime} = T_h$. 
We multiply (\ref{2.5a}) by $T_h(u_n)$
 and integrate by  parts to obtain
$$
\begin{array}{rl}
 \displaystyle\frac{d}{dt}   \int_{\Omega} \Phi_h(u_n)dx +
  \int_{ -h<u_n<0}  |\nabla u_n|^2dx =&\displaystyle
  \chi \int_{ -h<u_n<0} u_n \frac{|\nabla v_n |^{p-2}\nabla v_n \nabla u_{n}}{1+ \frac{1}{n} |\nabla v_n|^{p-1}}dx
  \smallskip\\
[4mm] 
 & \displaystyle
 +\mu \int_{\Omega} T_h(u_n)u_n(1-(u_n)_+)dx
 \\
[4mm] 
 \le  &\displaystyle
 \chi\, n  \int_{ -h<u_n<0} u_n  | \nabla u_{n} |dx
  +\mu \int_{\Omega} T_h(u_n)u_ndx
\end{array}
$$
and apply H\"older inequality and after some computations we get 
$$\frac{d}{dt}   \int_{\Omega} \Phi_h(u_n)dx +
  \frac{1}{2}  \int_{ -h<u_n<0}   |\nabla u_n|^2dx \leq 
  \frac{\chi^2}{2} n \int_{ -h<u_n<0}  u_n^2  
  +2\mu  \int_{\Omega}  \Phi_h(u_n) dx
  $$
  which implies 
  $$\frac{d}{dt}   \int_{\Omega} \Phi_h(u_n)dx +
   \frac{1}{2}\int_{ -h<u_n<0}   |\nabla u_n|^2dx \leq 
   ch^2  
  +2 \mu  \int_{\Omega}  \Phi_h(u_n)dx. 
  $$
  We divide by $h$ and take limits when $h \rightarrow 0$ to get 
    $$\frac{d}{dt}   \int_{\Omega} \ (-u_n)_+dx   \leq 
   2 \mu    \int_{ \Omega}  (-u_n)_+  dx .
  $$
  We apply Gronwall's Lemma and it results 
  $$\int_{\Omega} (-u_n)_+ dx \leq e^{2 \mu t}\int_{\Omega} (-u_n(0))_+ dx 
 .
$$
Since $ (-u_0)_+=0$,  we get 
$$ u_n \geq 0.$$  
Since $u_n \geq 0$ we have that $(u_n)_+= u_n$ and  $u_n$ satisfies 
 (\ref{2.1}). $\hfill\square$
\end{proof}\vspace{0.5cm}

The following estimate will be useful to obtain $L^q$ estimates for $u_n$. 

\begin{lemma}\label{l2.4}
Let $p<\frac{3}{2}$, $N\ge 3$ and $q$ satisfying assumptions 
 $$ \max\left\{  \frac{1}{2}, \frac{N-1}{2}- \frac{N}{4(p-1)} \right\} <q <\frac{N-1}{2},$$ then  
$$ \int_{\Omega} u_n^{2q}dx<c$$
and 
$$
\int_0^t \int_{\Omega} u_n^{2q+1} dx\,ds \leq  \overline{c}t+c_0.
 $$ 
\end{lemma}

\begin{proof}
We multiply equation (\ref{2.1}) by $|u_n|^{2q-2}u_n$  and after integration over $\Omega$ we get 
\begin{equation} \label{q2} \begin{array}{l} 
   \displaystyle 
   
   \frac{d}{dt}
   \frac{1}{2q}  \int_{\Omega} 
    u_n^{2q} dx+ \frac{2q-1}{q^2}\int_{\Omega} |\nabla u_n^q|^2 dx =  \frac{(2q-1)\chi }{q}\int_{\Omega} u_n^q  \frac{|\nabla v_n|^{p-2}\nabla v_n}{1+\frac{1}{n} |\nabla  v_n|^{p-1}}\nabla u_n^q dx + \mu \int_{\Omega}  u_n^{2q}(1-u_n) dx
    .\end{array} \end{equation} 
    Notice that 
    $$
    \frac{(2q-1)\chi}{q}
    \int_{\Omega} u_n^q  \frac{|\nabla v_n|^{p-2}\nabla v_n}{1+\frac{1}{n} |\nabla  v_n|^{p-1}}\nabla u_n^q dx
  \leq 
 \frac{(2q-1)}{2q^2} \int_{\Omega} |\nabla u_n^q|^2 dx
    +\frac{(2q-1) \chi^2}{2} \int_{\Omega} |\nabla v_n|^{2p-2}u_n^{2q} dx 
    $$
    and thanks to H\"older inequality 
    $$
    \begin{array}{ll}
   \displaystyle \int_{\Omega} u_n^{2q} |\nabla v_n|^{2p-2}dx & \displaystyle\leq 
    \left[ \int_{\Omega} |\nabla v_n |^{2(p-1)(2q+1)} dx \right]^{\frac{1}{2q+1}} \left[\int_{\Omega} 
    u_n^{(2q+1)}dx\right]^{\frac{2q}{2q+1}} 
    \smallskip\\
    &\displaystyle\leq \| \nabla v_n \|_{L^{2(p-1)(2q+1)}(\Omega)}^{2(p-1)} \|u_n \|^{2q}_{L^{2q+1}(\Omega)}.  
    \end{array}$$
    The elliptic regularity of the problem gives the inequality  
    $$\|\nabla v_n\|_{L^{2(p-1)(2q+1) }(\Omega)}
    \leq c \| u_n \|_{L^{2q+1}(\Omega)}
    $$ 
    provided
    $$ \frac{N (2q+1)}{N-2q-1}> 2(p-1)(2q+1)$$
  which is equivalent to  
   \begin{equation} \label{q1}   \frac{N}{N-2q-1}> 2(p-1) 
   \end{equation}
  and  satisfied for 
  $$\frac{N}{2(p-1)} > N-2q-1
  $$
  i.e. 
    $$\frac{N[2p-3]}{2(p-1)}<2q+1.
  $$
  We take $q$ such that 
  $$ \frac{N[2p-3]}{2(p-1)}< 2q+1$$
  which is  equivalent to 
 
      $$\frac{(N-1)[2p-2]-N}{2(p-1)} <2q
  $$
   \begin{equation}\label{2.5_0}
   \frac{N-1}{2} -\frac{N}{4(p-1)} <q.
   \end{equation}
  Thanks to \eqref{q1} and \eqref{2.5_0} we have that 
  $$q \in \left(
    \max\left\{  \frac{1}{2}, \frac{N-1}{2}- \frac{N}{4(p-1)} \right\} , \frac{N-1}{2}\right).$$
Then, after some computations, 
we get  
$$\frac{(2q-1) \chi}{q}  \int_{\Omega} u_n^q |\nabla v_n|^{p-2}\nabla v_n \nabla u_n^q dx\leq  \frac{(2q-1)}{2 q^2} \int_{\Omega} |\nabla u_n^q|^2 dx
    +c \left[  \int_{\Omega} u_n^{2q+1}dx
    \right]^{1+\frac{2(p-\frac{3}{2})}{2q+1}} 
    $$
    we replace the previous inequality  in equation (\ref{q2}) to get 
\begin{equation}  \label{eq2706}
    \frac{d}{dt} 
   \frac{1}{2q}  \int_{\Omega} 
    u_n^{2q} dx+ \frac{2q-1}{2q^2}\int_{\Omega} |\nabla u_n^q|^2 dx \leq c \left[  \int_{\Omega} u_n^{2q+1} dx
    \right]^{1+\frac{2(p-\frac{3}{2})}{2q+1}} -\frac{\mu}{2} \int_{\Omega} u_n^{2q+1}dx +c. 
\end{equation} 
We apply Poincar\'e-Wirtinger inequality to 
the term  $$ \int_{\Omega} | u_n^{2q}| dx $$
in the following way 
$$\int_{\Omega} | u_n^{2q}| dx \leq c\int_{\Omega} |\nabla u_n^q|^2 dx + c\left[\int_{\Omega} | u_n|^q dx\right]^2$$
then, 
$$\frac{2q-1}{2q^2}\int_{\Omega} |\nabla u_n^q|^2 dx
\geq c \int_{\Omega} | u_n^{2q}| dx- c\left[\int_{\Omega} | u_n|^q dx\right]^2.
$$
Since 
$$\left[\int_{\Omega} | u_n|^q dx\right]^2
\leq \epsilon \int_{\Omega} | u_n|^{2q+1} dx+ c(\epsilon)
$$
and 
$$\left[\int_{\Omega} | u_n|^{2q+1}   dx\right]^{1+ \frac{2(p- \frac{3}{2})}{2q+1}}
\leq \epsilon \int_{\Omega} | u_n|^{2q+1} dx+ c(\epsilon).
$$
Thanks to (\ref{eq2706}), we have 
         \begin{equation} \label{q4} \begin{array}{l} \displaystyle
    \frac{d}{dt} 
   \frac{1}{2q}  \int_{\Omega} 
    u_n^{2q} dx+      c \int_{\Omega} u_n^{2q} dx  \leq 
 - (\frac{\mu}{2}-2 \epsilon)   \int_{\Omega}  |u_n|^{2q+1} dx + c(\epsilon)
    .\end{array} \end{equation} 
  We take $\epsilon<\mu/4$ and  apply maximum principle for ODE to obtain
 $$  \int_{\Omega} u_n^{2q} dx \leq  c, 
 \qquad 
\int_0^t \int_{\Omega} u_n^{2q+1} dx\, ds \leq  \overline{c}t+c_0,
 $$  
 where $c$ and $\overline{c}$ are positive constants. 
 $\hfill\square$
  \end{proof}
  \vspace{0.5cm}
  
  In the result below we give an $L^{2q}$-estimate of $u_n$ for $N=2$.
 
  \begin{lemma}\label{lema2_5} Let 
  $N=2$,   $p<2$  and $q \in(1,\infty)$ then, 
  $$\int_{\Omega} u_n^{2q} dx <c(q).$$
  and 
  $v_n \in  W^{1, \infty}(\Omega).$
  \end{lemma}
  \begin{proof}
  We proceed as in Lemma \ref{l2.4} to obtain 
  \begin{equation} \label{q21} \begin{array}{l} 
   \displaystyle \frac{d}{dt}
   \frac{1}{2q}  \int_{\Omega} 
    u_n^{2q} dx+ \frac{2q-1}{2q^2}\int_{\Omega} |\nabla u_n^q|^2 dx \le 
   \displaystyle c\int_{\Omega} u_n^{2q} |\nabla v|^{2(p-1)} dx - \frac{\mu}{2} \int_{\Omega}  u_n^{2q+1} dx
   + c
    .\end{array} \end{equation} 
   Since 
  $$ \int_{\Omega} u_n^{2q} |\nabla v_n|^{2(p-1)} dx \leq 
    c \|\nabla v_n \|_{L^{\infty }(\Omega)}^{2(p-1)}  \int_{\Omega} 
    u_n^{2q}dx
    $$
and thanks to the  elliptic regularity we know that 
    $$\|\nabla v_n\|_{L^{ \infty}(\Omega)}
    \leq c \| u_n \|_{L^{2+\frac{\epsilon}{q}}(\Omega)}
    $$ 
    and for $q\geq 1$ 
    $$
    \| u_n \|_{L^{2+\frac{\epsilon}{q}}(\Omega)} \leq c \| u_n \|_{L^{2q+\epsilon }(\Omega)}
    $$ 
    then 
 \begin{equation} \label{eqwe}  \begin{array}{l} 
   \displaystyle \frac{d}{dt}
   \frac{1}{2q}  \int_{\Omega} 
    u_n^{2q} dx+ \frac{2q-1}{2q^2}\int_{\Omega} |\nabla u_n^q|^2 dx \leq 
   \displaystyle c \left[\int_{\Omega} u_n^{2q+\epsilon} dx \right]^{ 2\frac{p-1 +q }{2q+\epsilon}}
   - \frac{\mu}{2} \int_{\Omega}  u_n^{2q+1} dx
   + c
    .\end{array} \end{equation}
   Thanks to Gagliardo Nirenberg Inequality we know that 
 \begin{equation} \label{eq22}  \|u_n \|_{L^{2q+\epsilon}(\Omega)} =\|u_n^q \|_{L^{2+ \frac{\epsilon}{q}}(\Omega)}^{\frac{1}{q}}
   \leq c \left[ \| \nabla u_n^q \|_{L^2(\Omega) }^{\frac{a}{q}} \| u_n^q \|_{L^1(\Omega) }^{\frac{(1-a)}{q}}
   + c\| u_n^q \|_{L^1(\Omega) }\right]^{\frac{1}{q}}
   \end{equation}
   for $a$ satisfying 
   $$ \frac{1}{2+ \frac{\epsilon}{q}}= \left( \frac{1}{2}- \frac{1}{N} \right)a +(1-a)  $$
   which is equivalent to 
   $$a= \frac{2N(q+\epsilon)}{(2q+\epsilon)(N+2)}.$$
   For $N=2$ and $q=1$ we have that 
   $$a=\frac{1+\epsilon}{(2+\epsilon)}.$$
   We replace $q = 1 $ into (\ref{eq22}), and thanks to Lemma \ref{l3.1}, we get 
    $$ \|u_n \|_{L^{2+  \epsilon} (\Omega)}^{ 2+\epsilon }
   \leq 
  c  \| \nabla u_n \|_{L^2(\Omega) }^{1+\epsilon  }
  +c.
   $$
  Then, 
 $$ \left[\int_{\Omega} u_n^{2+\epsilon} dx \right]^{\frac{2p }{2+\epsilon}  }
  \leq c   \| \nabla u_n \|_{L^2(\Omega) }^{
  \frac{(1+ \epsilon) 2p }{2+ \epsilon}
  } +c,
  $$ 
  for $\epsilon$ small enough we get 
  $$ \| \nabla u_n \|_{
  L^2(\Omega)
  }^{
  \frac{(1+ \epsilon) 2p }{2+ \epsilon}
  }
  \leq \delta \| \nabla u_n\|_{L^2}^2 + c(\delta)
  $$
  for any $\delta>0.$
  We replace into (\ref{eqwe}) to get 
  $$ \begin{array}{l} 
   \displaystyle \frac{d}{dt}
   \frac{1}{2}  \int_{\Omega} 
    u_n^{2} dx+ \left( \frac{1}{2} - \delta\right) \int_{\Omega} |\nabla u_n|^2 dx~ \le 
   \displaystyle c - \frac{\mu}{4} \int_{\Omega}  |u_n|^{3} dx
    .\end{array} 
    $$
  Standard computations shows 
  the uniform boundedness in time  of $$\int_{\Omega} 
    u_n^{2} dx.$$ 
    To end the proof we take $q\in(1,2)$ to obtain,  thanks to Gagliardo-Nirenberg inequality that 
    \begin{equation} \label{eq2222}  \|u_n \|_{L^{2q+\epsilon}(\Omega)} =\|u_n^q \|_{L^{2+ \frac{\epsilon}{q}}(\Omega)}^{\frac{1}{q}}
   \leq c \left[ \| \nabla u_n^q \|_{L^2(\Omega) }^{a} \| u_n^q \|_{L^{\frac{2}{q}}(\Omega) }^{1-a}
   + c\| u_n^q \|_{L^1(\Omega) }\right]^{\frac{1}{q}}
   \end{equation}
   for $a$ satisfying 
   $$ \frac{1}{2+ \frac{\epsilon}{q}}= \left( \frac{1}{2}- \frac{1}{N} \right)a +\frac{q(1-a)}{2}  $$
   i.e. 
    $$ \frac{q}{2q+ \epsilon}-\frac{q}{2}= \left( \frac{1}{2}- \frac{1}{N} \right)a -\frac{qa}{2}  $$
    since $N=2$ we have 
     $$ \frac{q}{2q+ \epsilon}-\frac{q}{2}=  -\frac{qa}{2}  $$
     and then
       $$ \frac{1}{2q+ \epsilon}-\frac{1}{2}=  -\frac{a}{2}  $$
         $$ \frac{2}{2q+ \epsilon}-1=  -a  $$
    i.e 
      $$ \left(  \frac{2(q-1)+ \epsilon}{2q+ \epsilon} \right)= a   $$
      then, for $\epsilon$ small enough and $q$ close to $1$, and since $p<2$, we have that 
      $$ a \left( 2\frac{p-1+q}{2q+ \epsilon}\right) <2.$$ 
      We replace into (\ref{eqwe}) to obtain 
      $$  \begin{array}{l} 
   \displaystyle \frac{d}{dt}
   \frac{1}{2q}  \int_{\Omega} 
    u_n^{2q} dx+ \frac{2q-1}{2q^2}\int_{\Omega} |\nabla u_n^q|^2 dx \leq 
   c \left[\int_{\Omega} |\nabla u_n^{q}|^2   dx \right]^{  \beta }
   - \frac{\mu}{2} \int_{\Omega}  u_n^{2q+1} dx
   + c
    \end{array} $$
    for $\beta<2$. We proceed as before to  get  $$\int_{\Omega} u_n^{2q} dx<c$$
    for some $q>1$. We apply elliptic regularity to second equation and get that $v_n \in W^{1, \infty}(\Omega)$. 
    Then, the term 
    $$\int_{\Omega} u_n^{2q}|\nabla v_n|^{2(p-1)} \leq c  \int_{\Omega} u_n^{2q}\leq c(\epsilon )+ \epsilon  \int_{\Omega} u_n^{2q+1}$$
    which implies, for 
    for $\epsilon \leq \frac{\mu}{4}$
 $$ \begin{array}{l} 
   \displaystyle \frac{d}{dt}
   \frac{1}{2q}  \int_{\Omega} 
    u_n^{2q} dx+ \frac{2q-1}{2q^2}\int_{\Omega} |\nabla u_n^q|^2 dx \leq 
   c - \frac{\mu}{4} \int_{\Omega}  u_n^{2q+1} dx
    .\end{array} 
    $$
    Since 
    $$\int_{\Omega} 
u^{2q}dx \leq 
\int_{\Omega} 
u^{2q+1}dx+ |\Omega|$$ we get 
$$
  \frac{d}{dt}
   \frac{1}{2q}  \int_{\Omega} 
    u_n^{2q} dx+ \frac{2q-1}{2q^2}\int_{\Omega} |\nabla u_n^q|^2 dx
    + \frac{\mu}{4}
    \int_{\Omega}  u_n^{2q} dx
    \leq 
   c +|\Omega|.
    $$
Maximum principle ends the proof. $\hfill\square$
  \end{proof}
  \vspace{0.5cm}

In the following Lemma   we obtain an estimate of $u_n$ and $v_n$, for $N\ge 3$. 
    \begin{lemma} \label{l4.10}
  Let  $N \geq 3$ and  $p<\frac{3}{2}$, then, for any  $s<\infty$   we have that
  $$ \int_{\Omega} u_n^{s} \leq c$$ and 
  $v_n \in W^{2, N+1} (\Omega) \cap W^{1, \infty}(\Omega).$
  \end{lemma}
  \begin{proof} Thanks to Lemma \ref{l2.4}, we have that $u_n$ is uniformly bounded in  
  $ L^{q} (\Omega)$ for $q<N-1$, then  we have that, thanks to the elliptic regularity,    
   $v_n\in W^{2,q} (\Omega) $ for $q<N-1$,
  which implies $$v_n \in  W^{1,\frac{Nq}{N-q}} (\Omega)$$
 for any $q\in (1, N-1)$.  In particular,  we have that, for any  $p<3/2$, $q<N-1$ and close to $N-1$, the inequality  
 $$ 2(p-1)(2q+1) <   \frac{Nq}{N-q}$$
 is satisfied and then,  the term  
 $$\|\nabla v_n\|_{L^{ 2(p-1)(2q+1) }(\Omega)}$$ is bounded. Then we proceed as in the case $N=2$ (see Lemma \ref{lema2_5})  to get   that 
 $u_n$ is uniformly bounded in $L^{s} (\Omega)$ for any $s<\frac{N(N-1)}{2(p-1)}-1$. Since $p<\frac{3}{2}$ we get that $u_n$ is uniformly bounded in $L^{s} (\Omega)$ for any $s <  N(N-1) -1$ and we deduce,  in view of $N \geq 3$ that $u_n \in L^{N+1}(\Omega)$. Elliptic regularity implies that $v_n \in W^{2, N+1} (\Omega)$ which is included in $W^{1, \infty}(\Omega)$  and the proof ends. 
 $\hfill\square$\end{proof}\vspace{0.5cm}

 The following result gives a uniform bound of $u_n$ when $N\ge 3$ and $p<3/2$.
 \begin{lemma}\label{lema_2.8}
  Let  $N \geq 2$ and  $p<\frac{3}{2}$, then,    we have that 
  $$ \| u_n \|_{L^{\infty}(\Omega)} \leq c.$$
  \end{lemma}
  \begin{proof}
 Thanks to Lemma \ref{lema2_5} and \ref{l4.10}, we have that   $v_n\in W^{1, \infty}(\Omega)$. We take $u_n^{2q-1}$ as test function in the weak formulation of (\ref{2.1})-(\ref{2.4}) to get, after some computations, that   
\begin{equation} \label{nueva} 
  \frac{1}{2q} \frac{d}{dt} 
  \int_{\Omega} u_n^{2q}dx + \frac{2q-1}{2q^2} 
  \int_{\Omega} |\nabla u_n^{q}|^2 dx\leq \|\nabla v_n\|_{L^{\infty} (\Omega) }^{2(p-1)}\frac{\chi^2(2q-1)}{2} \int_{\Omega} u_n^{2q}dx+ \mu \int_{\Omega} u_n^{2q}dx-\mu \int_{\Omega} u_n^{2q+1}dx.
   \end{equation} 
 We split  the term 
 $ \mu \int_{\Omega} u_n^{2q}dx$ into two parts as follows   
  $$ \mu \int_{\Omega} u_n^{2q}dx
 \leq 
 \frac{\mu}{2} \int_{\Omega} u_n^{2q+1}dx+
 \mu|\Omega |4^q.
$$
In the case $N=2$, we apply Gagliardo-Nirenberg inequality to the term  
  $$
  \int_{\Omega} u_n^{2q}dx
  $$
in the following way 
$$  \|u_n \|_{L^{2q}(\Omega)} =\|u_n^q \|_{L^2(\Omega)}^{\frac{1}{q}}
   \leq c \left[ \| \nabla u_n^q \|_{L^2(\Omega) }^{\frac{a}{q}} \| u_n^q \|_{L^1(\Omega) }^{\frac{(1-a)}{q}}
   + c\| u_n^q \|_{L^1(\Omega) }\right]^{\frac{1}{q}}
 $$
   for $a$  defined by 
   $$a= \frac{N}{(N+2)}.$$
Then, thanks to Young$^{\prime}$s inequality 
  $$   \int_{\Omega} u_n^{2q}dx\leq 
  \epsilon  \int_{\Omega} u_n^{2q+1}dx
  +   c\epsilon^{-\frac{N}{2}} 
  $$
  for any $\epsilon>0$. We take
  $$
  \epsilon < \frac{2q-1}{4q^2} 
  \frac{2}{ \|\nabla v_n\|_{L^{\infty} (\Omega) }^{2(p-1)}\chi^2(2q-1)} 
  =
  \frac{1}{4q^2} 
  \frac{2}{ \|\nabla v_n\|_{L^{\infty} (\Omega) }^{2(p-1)}\chi^2}
  $$ 
  we have that 
  $$
  \|\nabla v_n\|_{L^{\infty} (\Omega) }^{2(p-1)}\frac{\chi^2(2q-1)}{2}    \int_{\Omega} u_n^{2q}dx
  \leq 
  \frac{2q-1}{4q^2}  \int_{\Omega} |\nabla u_n^{q}|^2dx
  +  c(1+ q^N).
  $$
  We replace the previous inequality into 
  (\ref{nueva}) to get, in view of 
  $$u_n^{2q} \leq u_n^{2q+1} +1$$
  and 
  $$c(2q-1)u_n^{q} \leq \frac{\mu}{4} u_n^{2q} +cq^2.$$
  For the case $N \geq 3$, we have that 
 $$
   \frac{1}{2q} \frac{d}{dt} 
  \int_{\Omega} u_n^{2q}dx + \frac{2q-1}{4q^2} 
  \int_{\Omega} |\nabla u_n^{q}|^2 dx
  +\frac{\mu}{2} \int_{\Omega} u_n^{2q} dx
  \leq  c (q^{N} +1).
 $$ Thanks to maximum principle we have that 
 $$\int_{\Omega} u_n^{2q} dx
  \leq  c (q^{N}+1). 
 $$
    We take q-roots in the previous expression to get 
$$  \|u_n\|_{L^{2q}(\Omega)} \leq  c$$

for some  $c$  independent of $q$. We take limits when $q \rightarrow \infty$ to end the proof. 
$\hfill\square$
\end{proof}\vspace{0.5cm}

 \begin{lemma} 
 Let $u_n$ be the solution to (\ref{2.1})-(\ref{2.4}) satisfying assumptions 
 (\ref{H1})-(\ref{H3}), 
 then, for any $T \in(0, \infty)$, we have  that $u_n$ satisfies 
 $$\int_0^T  
 \int_{\Omega} |\nabla u_n|^2 dxds\leq c_1T +c_2. $$
 \end{lemma} 
 \begin{proof} 
 We denote by $y$ the average of $u_n$, i.e. 
 $$y(t):= \frac{1}{|\Omega|}  \int_{\Omega} u_n dx$$
 which satisfies 
 $$y^{\prime}=\mu \int_{\Omega} u_n(1-u_n) dx.$$  
We  multiply equation (\ref{2.1}) by $(u_n-y)$ and integrate by parts to get,  
 $$
\begin{array}{l} 
   \displaystyle \frac{d}{dt} \int_{\Omega} | u_n-y|^2 dx
   +\int_{\Omega} |\nabla u_n|^2 dx =   \chi \int_{\Omega} u_n \frac{|\nabla v_n|^{p-2}\nabla v_n}{1+\frac{1}{n} |\nabla  v_n|^{p-1}} \nabla u_n\, dx  + \mu \int_{\Omega}  u_n(1-u_n)( u_n-y)dx
    .\end{array} 
$$    
Thanks to Lemma \ref{lema_2.8}
$$ 
\int_{\Omega}  u_n(1-u_n)(u_n-y)dx 
\leq c
$$
and 
$$
\begin{array}{ll}
\displaystyle\int_{\Omega} u_n \frac{|\nabla v_n|^{p-2}\nabla v_n}{1+\frac{1}{n} |\nabla  v_n|^{p-1}} \nabla u_n dx &\displaystyle \leq \frac{1}{2} \int_{\Omega} |\nabla u_n|^2 dx
    +c \int_{\Omega} u_n^{2} |\nabla v_n|^{2p-2}dx $$
   \smallskip\\
   &\displaystyle\leq \frac{1}{2} \int_{\Omega} |\nabla u_n|^2 dx+
    c .
    \end{array}
    $$
After integration we get the wished result. $\hfill\square$
  \end{proof}
  \vspace{0.5cm}

\begin{lemma}\label{lemma2.9} Let $u_n$ the solution to (\ref{1.1})-(\ref{1.4}) then, under assumptions (\ref{H1})-(\ref{H3}) and for any   for any $T \in(0, \infty)$,    we have  that
$$u_{nt} \in L^2(0,T: (H^1(\Omega))^{\prime}).$$  
\end{lemma}
\begin{proof} For simplicity, we denote by $X$ the space $L^2(0,T: H^1(\Omega))$ and by $X^{\prime} $ its dual, which is equivalent to 
$L^2(0,T: (H^1(\Omega))^{\prime})$. 

Since  $$u_{nt} = \Delta u_n - div \left(\chi u_n \frac{|\nabla v_n |^{p-1}\nabla v_n}{1+\frac{1}{n}|\nabla v_n|}\right)  +
u_n(1-u_n) $$
and for any $w \in L^2(0,T:H^1(\Omega))$  we have that 
\begin{equation} \label{julio1} \left|<w, -\Delta u_n>_{X, X^{\prime}} \right|=\left|\int_0^T\int_{\Omega} 
 \nabla  u_n \nabla w dx dt \right| \leq \|u_n \|_{L^2(0,T:H^1(\Omega))}  
 \|w \|_{L^2(0,T:H^1(\Omega))}  ,
 \end{equation} 
 $$
 \begin{array}{ll}
 \left|<w, div \left(\chi u_n \frac{|\nabla v_n |^{p-1}\nabla v_n}{1+\frac{1}{n}|\nabla v_n|}\right)>_{X, X^{\prime}} \right|&\displaystyle=\left|\int_0^T\int_{\Omega} 
\chi u_n \frac{|\nabla v_n |^{p-1}\nabla v_n}{1+\frac{1}{n}|\nabla v_n|} \nabla w dx dt \right|
\smallskip\\
&\displaystyle\leq c  
 \|w \|_{L^2(0,T:H^1(\Omega))}
 \end{array}
 $$
 which implies 
 \begin{equation} 
 \label{julio2}
 \left|<w, div \left(\chi u_n \frac{|\nabla v_n |^{p-1}
 \nabla v_n}{1+\frac{1}{n}|\nabla v_n|}\right)>_{X, X^{\prime}} \right| 
\leq c  
 \|w \|_{L^2(0,T:H^1(\Omega))}.
 \end{equation}
Finally,  
$$ 
\begin{array}{ll}
\left|< w, u_n(1-u_n)>_{X, X^{\prime}}\right| & =
 \displaystyle\left|\int_0^T \int_{\Omega}
u_n(1-u_n) wdx  \right|
\smallskip\\
&\leq c \|w \|_{L^2(0,T:H^1(\Omega))} 
\end{array}
$$
which gives 
\begin{equation}
    \label{julio3}
\left|< w, u_n(1-u_n)>_{X, X^{\prime}}\right|  \leq c \|w \|_{L^2(0,T:H^1(\Omega))} 
\end{equation} 
From (\ref{julio1}), (\ref{julio2}) and (\ref{julio3})   we deduce the result.
$\hfill\square$
  \end{proof}
  \vspace{0.5cm}

 \section{Existence of solutions} 
 \label{s3}
 \setcounter{equation}{0}

The proof of the existence of solutions is given into several steps. First we prove the existence of weak solutions of the approximated problem (\ref{2.1})-(\ref{2.4}) and thanks to the estimates obtained in  Section \ref{s2}, we prove the convergence    of the weak formulation of 
 (\ref{2.1})-(\ref{2.4}) to the weak  solution of (\ref{1.1})-(\ref{1.4}).   
  
 \begin{lemma}
     Let  $p<3/2$, then, under assumptions (\ref{H1})-(\ref{H3}) 
      there exists a solution for the weak formulation of the approximated problem (\ref{2.1})-(\ref{2.4}) satisfying 
      $$u_n \in L^2(0,T: L^2(\Omega)). 
      $$
 \end{lemma}
\begin{proof} We consider a  fixed point argument and 
define the following set 
$$ A:= \{ w \in L^2(0,T:L^2(\Omega)), \quad w\geq 0 \quad  \|w \|_{L^{\infty}(0,T:L^{\infty}(\Omega))} 
\leq c(T) \} 
$$
where $c(T)$ is the constant obtained in Lemma \ref{lema2_5} for $N=2$ and Lemma \ref{l4.10} for $N \geq 3$. Now, we  consider the function $J$ 
$$J: A \subset  L^2(0,T: L^2(\Omega)) \rightarrow L^2(0,T:L^2(\Omega)).$$
  Let $\tilde{u}_{n} \in L^2(0,T: L^2(\Omega))$ and define $J(\tilde{u}_{n})= u_{n} $ as the solution to the problem 
$$
{u_n}_t-\Delta u_{n}  = -  div  \left(\chi u_{n}|\frac{\nabla v_{n}|^{p-2}\nabla v_{n}}{1+\frac{1}{n}|\nabla v_n|}\right)+\mu u_{n}( 1- \tilde{u}_{n})
$$
where $v_n$ is the solution to 
$$-\Delta v_n+ v_n  = \tilde{u}_{n}
.$$
We first notice that $J$ is a continuous function and 
thanks to the estimates obtained in Section \ref{s2}, we have that  
\begin{itemize} 
\item[(i)] 
$ J(A)$ is a precompact set in $A$;
\item[(ii)] 
$J(A) \subset A$; 
\item[(iii)]
$J(A)$ is a bounded set in $H^1(0,T: (H^1(\Omega))^{\prime}) \cap L^2(0,T:H^1(\Omega))$.
\end{itemize} 
Thanks to Aubin Lions' Lemma we get the result for any $T<\infty$. 
$\hfill\square$
  \end{proof}
  \vspace{0.5cm}
  
 \begin{lemma}  Let  $p<3/2$, then, under assumptions (\ref{H1})-(\ref{H3}) the weak solution of (\ref{2.1})-(\ref{2.4}) converges to the weak solution of 
 (\ref{1.1})-(\ref{1.4}).     %
 \end{lemma}
\begin{proof}
   We reproduce the steps given in   section \ref{s2}, to   obtain  that
    \[
    \begin{array}{lll}
        \displaystyle\int_\Omega |\nabla u_{n}|^2 dx&\le C,&
        \smallskip\\
     \displaystyle\|u_{n}\|_{L^{\infty}(\Omega)}
        &\leq C, &\displaystyle
        \smallskip\\
        \| v_{n}\|_{W^{2, q}(\Omega)} 
        &\leq C &   \mbox{ for any}~ q<\infty,  \smallskip\\ 
     \|u_{n t}\|_{L^{2}(0,T:(H^1(\Omega))^{\prime})}
        &\leq C(T+1). &
    \end{array}
    \]
    Since $u_n \in  H^{1}(0,T:(H^1(\Omega))^{\prime})$ we have the 
    $v_n \in  H^{1}(0,T:(H^1(\Omega)))$.  We consider the inclusions of spaces $$W^{2,N+1}(\Omega) \subset W^{1, \infty}(\Omega) \subset H^1( \Omega).$$ Since   
    $$ v_n \in L^{\infty}(0,T:W^{2, N+1}) \cap 
    H^1(0,T:H^{1}(\Omega)),
$$ we apply Aubin-Lions Lemma to 
    $v_n$ and deduce that 
    $$v_n \in C(0,T: W^{1, \infty}(\Omega)$$
    and there exists $(u^*,v^*)$ and a  subsequence $(u_{n_j},  v_{n_j}) $ such that 
    $$v_{n_j} \longrightarrow v^* , \mbox{ strongly   in } 
   L^p(0,T:W^{1, \infty}(\Omega)), \quad \mbox{ for any } p<\infty,$$
which implies that 
$$ \frac{1}{1+ \frac{1}{n_j} |\nabla v_{n_j} |}   \longrightarrow 1, 
\mbox{ weakly    in }  L^4(0,T:L^4(\Omega)); $$
$$  |\nabla v_{n_j}|^{p-2} \nabla v_{n_j}   \longrightarrow |\nabla v^{*}|^{p-2} \nabla v^*, 
\mbox{ strongly     in }  L^4(0,T:L^4(\Omega)). $$
Moreover
    \[
    \begin{array}{ll}
   u_{n j} \longrightarrow u^* , &
   \mbox{ strong  in } L^2(0,T:L^2(\Omega)).
    \\
    u_{n_j} \rightharpoonup u^* , & \mbox{weak  in } L^2(0,T:H^1(\Omega))
    . 
     \end{array}
     \]
    We take limits in the weak formulation of the approximated problem in Definition \ref{def3.1} to obtain that $(u^*, v^*)$ satisfies Definition \ref{def1.1} and proves the existence of weak solutions.
\end{proof}

\end{document}